\documentclass[11pt]{article}
\usepackage{amssymb,amsthm,amsmath,amscd}

\newtheorem{thm}{Theorem}[section]
\newtheorem{prop}[thm]{Proposition}
\newtheorem{lem}[thm]{Lemma}
\newtheorem{cor}[thm]{Corollary}

\theoremstyle{definition}
\newtheorem{df}[thm]{Definition}
\newtheorem{rem}[thm]{Remark}

\renewcommand{\phi}{\varphi}

\numberwithin{equation}{section}

\newcommand{\N}{\mathbb{N}}
\newcommand{\Z}{\mathbb{Z}}
\newcommand{\Q}{\mathbb{Q}}
\newcommand{\R}{\mathbb{R}}
\newcommand{\C}{\mathbb{C}}
\newcommand{\T}{\mathbb{T}}

\newcommand{\Hom}{\operatorname{Hom}}
\newcommand{\Aut}{\operatorname{Aut}}
\newcommand{\aInn}{\operatorname{\overline{Inn}}}

\newcommand{\Ad}{\operatorname{Ad}}
\newcommand{\id}{\operatorname{id}}
\newcommand{\Ima}{\operatorname{Im}}
\newcommand{\diag}{\operatorname{diag}}
\newcommand{\Sp}{\operatorname{Sp}}
\newcommand{\Lip}{\operatorname{Lip}}
\newcommand{\OrderExt}{\operatorname{OrderExt}}
\newcommand{\Ext}{\operatorname{Ext}}
\newcommand{\Aff}{\operatorname{Aff}}

\newcommand{\ep}{\varepsilon}

\title{$\Z$-actions on AH algebras and $\Z^2$-actions on AF algebras}
\author{Hiroki Matui \\
Graduate School of Science \\
Chiba University \\
Inage-ku, Chiba 263-8522, Japan}
\date{}

\begin{document}
\maketitle

\begin{abstract}
We consider $\Z$-actions (single automorphisms) 
on a unital simple AH algebra with real rank zero and slow dimension growth 
and show that the uniform outerness implies the Rohlin property 
under some technical assumptions. 
Moreover, two $\Z$-actions with the Rohlin property 
on such a $C^*$-algebra are shown to be cocycle conjugate 
if they are asymptotically unitarily equivalent. 
We also prove that 
locally approximately inner and uniformly outer $\Z^2$-actions 
on a unital simple AF algebra with a unique trace 
have the Rohlin property 
and classify them up to cocycle conjugacy 
employing the $\OrderExt$ group as classification invariants. 
\end{abstract}

\section{Introduction}

Classification of group actions is 
one of the most fundamental subjects in the theory of operator algebras. 
For AFD factors, 
a complete classification is known for actions of countable amenable groups. 
However, classification of automorphisms or group actions 
on $C^*$-algebras is still a far less developed subject, 
partly because of $K$-theoretical difficulties. 
For AF and AT algebras, 
A. Kishimoto \cite{K95crelle,K96JFA,K98JOT,K98JFA} showed 
the Rohlin property for a certain class of automorphisms and 
obtained a cocycle conjugacy result. 
Following the strategy developed by Kishimoto, 
H. Nakamura \cite{N2} showed that 
aperiodic automorphisms on unital Kirchberg algebras are 
completely classified by their $KK$-classes 
up to $KK$-trivial cocycle conjugacy. 
As for $\Z^N$-actions, 
Nakamura \cite{N1} introduced the notion of the Rohlin property and 
classified product type actions of $\Z^2$ on UHF algebras. 
T. Katsura and the author \cite{KM} gave a complete classification 
of uniformly outer $\Z^2$-actions on UHF algebras 
by using the Rohlin property. 
For Kirchberg algebras, 
M. Izumi and the author \cite{IM} classified a large class of $\Z^2$-actions 
and also showed the uniqueness of $\Z^N$-actions 
on the Cuntz algebras $\mathcal{O}_2$ and $\mathcal{O}_\infty$. 
The present article is a continuation of these works. 

In the first half of this paper, 
we study single automorphisms (i.e. $\Z$-actions) 
on a unital simple classifiable AH algebra. 
More precisely, we prove the Rohlin type theorem (Theorem \ref{ZRohlintype}) 
for an automorphism $\alpha$ of a unital simple AH algebra $A$ 
with real rank zero and slow dimension growth 
under the assumption that $\alpha^r$ is approximately inner for some $r\in\N$ 
and $A$ has finitely many extremal tracial states. 
Furthermore, we also prove that 
if two automorphisms $\alpha$ and $\beta$ 
of a unital simple AH algebra with real rank zero and slow dimension growth 
have the Rohlin property and 
$\alpha\circ\beta^{-1}$ is asymptotically inner, then 
the two $\Z$-actions generated by $\alpha$ and $\beta$ are cocycle conjugate 
(Theorem \ref{Zcc}). 
These results are generalizations of 
Kishimoto's work for AF and AT algebras 
(\cite{K95crelle,K96JFA,K98JOT,K98JFA}). 
For the proofs, 
we need to improve some of the arguments in \cite{K98JOT,K98JFA} 
concerning projections and unitaries in central sequence algebras. 
As a byproduct, it will be also shown that for $A$ as above, 
the central sequence algebra $A_\omega$ satisfies 
Blackadar's second fundamental comparability question 
(Proposition \ref{byproduct}). 

In the latter half of the paper, 
we study $\Z^2$-actions on a unital simple AF algebra $A$ with a unique trace. 
We first show 
a $\Z$-equivariant version of the Rohlin type theorem for single automorphisms 
(Theorem \ref{Z2Rohlin5}), 
and as its corollary 
we obtain the Rohlin type theorem 
for a $\Z^2$-action $\phi:\Z^2\curvearrowright A$ 
such that $\phi_{(r,0)}$ and $\phi_{(0,s)}$ are approximately inner 
for some $r,s\in\N$ (Corollary \ref{Z2Rohlin6}). 
This is a generalization of \cite[Theorem 3]{N1}. 
Next, by using a $\Z$-equivariant version of 
the Evans-Kishimoto intertwining argument \cite{EK}, 
we classify uniformly outer locally $KK$-trivial $\Z^2$-actions on $A$ 
up to $KK$-trivial cocycle conjugacy (Theorem \ref{Z2cc}). 
This is a generalization of \cite[Theorem 6.5]{KM}. 
We remark that 
$KK$-triviality of $\alpha\in\Aut(A)$ is equivalent to 
$\alpha$ being approximately inner, 
because $A$ is assumed to be AF. 
For classification invariants, 
we employ the $\OrderExt$ group introduced in \cite{KK}. 
The crossed product of $A$ by the first generator $\phi_{(1,0)}$ is 
known to be a unital simple AT algebra with real rank zero. 
The second generator $\phi_{(0,1)}$ naturally extends to 
an automorphism of the crossed product and 
its $\OrderExt$ class gives the invariant of the $\Z^2$-action $\phi$. 
However we do not know the precise range of the invariant in general.

\section{Preliminaries}

We collect notations and terminologies relevant to this paper. 
For a Lipschitz continuous map $f$ between metric spaces, 
$\Lip(f)$ denotes the Lipschitz constant of $f$. 
Let $A$ be a unital $C^*$-algebra. 
For $a,b\in A$, we mean by $[a,b]$ the commutator $ab-ba$. 
For a unitary $u\in A$, 
the inner automorphism induced by $u$ is written by $\Ad u$. 
An automorphism $\alpha\in\Aut(A)$ is called outer, 
when it is not inner. 
A single automorphism $\alpha$ is often identified 
with the $\Z$-action induced by $\alpha$. 
An automorphism $\alpha$ of $A$ is said to be asymptotically inner, 
if there exists a continuous family of unitaries 
$\{u_t\}_{t\in[0,\infty)}$ in $A$ such that 
\[
\alpha(a)=\lim_{t\to\infty}\Ad u_t(a)
\]
for all $a\in A$. 
When there exists a sequence of unitaries $\{u_n\}_{n\in\N}$ in $A$ 
such that 
\[
\alpha(a)=\lim_{n\to\infty}\Ad u_n(a)
\]
for all $a\in A$, 
$\alpha$ is said to be approximately inner. 
The set of approximately inner automorphisms of $A$ is 
denoted by $\aInn(A)$. 
Two automorphisms $\alpha$ and $\beta$ are said to be 
asymptotically (resp. approximately) unitarily equivalent 
if $\alpha\circ\beta^{-1}$ is asymptotically (resp. approximately) inner. 
The set of tracial states on $A$ is denoted by $T(A)$. 
We mean by $\Aff(T(A))$ 
the space of $\R$-valued affine continuous functions on $T(A)$. 
The dimension map $D_A:K_0(A)\to\Aff(T(A))$ is 
defined by $D_A([p])(\tau)=\tau(p)$. 
For a homomorphism $\rho$ between $C^*$-algebras, 
$K_0(\rho)$ and $K_1(\rho)$ mean the induced homomorphisms on $K$-groups. 

As for group actions on $C^*$-algebras, 
we freely use the terminology and notation 
described in \cite[Definition 2.1]{IM}. 

Let $A$ be a separable $C^*$-algebra and 
let $\omega\in\beta\N\setminus\N$ be a free ultrafilter. 
We set 
\[
c_0(A)=\{(a_n)\in\ell^\infty(\N,A)\mid
\lim_{n\to\infty}\lVert a_n\rVert=0\},\quad 
A^\infty=\ell^\infty(\N,A)/c_0(A),
\]
\[
c_\omega(A)=\{(a_n)\in\ell^\infty(\N,A)\mid
\lim_{n\to\omega}\lVert a_n\rVert=0\},\quad 
A^\omega=\ell^\infty(\N,A)/c_\omega(A). 
\]
We identify $A$ with the $C^*$-subalgebra of $A^\infty$ (resp. $A^\omega$) 
consisting of equivalence classes of constant sequences. 
We let 
\[
A_\infty=A^\infty\cap A',\quad A_\omega=A^\omega\cap A'
\]
and call them the central sequence algebras of $A$. 
A sequence $(x_n)_n\in\ell^\infty(\N,A)$ is called a central sequence 
if $\lVert[a,x_n]\rVert\to0$ as $n\to\infty$ for all $a\in A$. 
A central sequence is a representative of an element in $A_\infty$. 
An $\omega$-central sequence is defined in a similar way. 
When $\alpha$ is an automorphism on $A$ or 
an action of a discrete group on $A$, 
we can consider its natural extension 
on $A^\infty$, $A^\omega$, $A_\infty$ and $A_\omega$. 
We denote it by the same symbol $\alpha$. 

Next, we would like to recall 
the definition of uniform outerness introduced by Kishimoto and 
the definition of the Rohlin property of $\Z^N$-actions 
introduced by Nakamura. 

\begin{df}[{\cite[Definition 1.2]{K95crelle}}]
An automorphism $\alpha$ of a unital $C^*$-algebra $A$ is 
said to be uniformly outer 
if for any $a\in A$, any non-zero projection $p\in A$ and any $\ep>0$, 
there exist projections $p_1,p_2,\dots,p_n$ in $A$ such that 
$p=\sum p_i$ and $\lVert p_ia\alpha(p_i)\rVert<\ep$ 
for all $i=1,2,\dots,n$. 
\end{df}

We say that 
an action $\alpha$ of a discrete group on $A$ is uniformly outer 
if $\alpha_g$ is uniformly outer for every element $g$ of the group 
other than the identity element. 

Let $N$ be a natural number. 
Let $\xi_1,\xi_2,\dots,\xi_N$ be the canonical basis of $\Z^N$, 
that is, 
\[
\xi_i=(0,0,\dots,1,\dots,0,0), 
\]
where $1$ is in the $i$-th component. 
For $m=(m_1,m_2,\dots,m_N)$ and $n=(n_1,n_2,\dots,n_N)$ in $\Z^N$, 
$m\leq n$ means $m_i\leq n_i$ for all $i=1,2,\dots,N$. 
For $m=(m_1,m_2,\dots,m_N)\in\N^N$, we let 
\[
m\Z^N=\{(m_1n_1,m_2n_2,\dots,m_Nn_N)\in\Z^N
\mid (n_1,n_2,\dots,n_N)\in\Z^N\}. 
\]
For simplicity, we denote $\Z^N/m\Z^N$ by $\Z_m$. 
Moreover, we may identify $\Z_m=\Z^N/m\Z^N$ with 
\[
\{(n_1,n_2,\dots,n_N)\in\Z^N
\mid 0\leq n_i\leq m_i{-}1\quad\forall i=1,2,\dots,N\}. 
\]

\begin{df}[{\cite[Definition 1]{N1}}]
Let $\phi$ be an action of $\Z^N$ on a unital $C^*$-algebra $A$. 
Then $\phi$ is said to have the Rohlin property, 
if for any $m\in\N$ there exist $R\in\N$ and 
$m^{(1)},m^{(2)},\dots,m^{(R)}\in\N^N$ 
with $m^{(1)},\dots,m^{(R)}\geq(m,m,\dots,m)$ 
satisfying the following: 
For any finite subset $F$ of $A$ and $\ep>0$, 
there exists a family of projections 
\[
e^{(r)}_g\qquad (r=1,2,\dots,R, \ g\in\Z_{m^{(r)}})
\]
in $A$ such that 
\[
\sum_{r=1}^R\sum_{g\in\Z_{m^{(r)}}}e^{(r)}_g=1,\quad 
\lVert[a,e^{(r)}_g]\rVert<\ep,\quad 
\lVert\phi_{\xi_i}(e^{(r)}_g)-e^{(r)}_{g+\xi_i}\rVert<\ep
\]
for any $a\in F$, $r=1,2,\dots,R$, $i=1,2,\dots,N$ 
and $g\in\Z_{m^{(r)}}$, 
where $g+\xi_i$ is understood modulo $m^{(r)}\Z^N$. 
\end{df}

It is clear that 
if $\phi:\Z^N\curvearrowright A$ has the Rohlin property, 
then $\phi$ is uniformly outer. 

We recall the definition of tracial rank zero introduced by H. Lin. 

\begin{df}[{\cite{LTR1,LTR2}}]
A unital simple $C^*$-algebra $A$ is said to have tracial rank zero 
if for any finite subset $F\subset A$, any $\ep>0$ and 
any non-zero positive element $x\in A$ 
there exists a finite dimensional subalgebra $B\subset A$ with $p=1_B$ 
satisfying the following. 
\begin{enumerate}
\item $\lVert[a,p]\rVert<\ep$ for all $a\in F$. 
\item The distance from $pap$ to $B$ is less than $\ep$ for all $a\in F$. 
\item $1-p$ is Murray-von Neumann equivalent to a projection in $xAx$. 
\end{enumerate}
\end{df}

In \cite{LDuke}, H. Lin gave a classification theorem 
for unital separable simple nuclear $C^*$-algebras 
with tracial rank zero which satisfy the UCT. 
(\cite[Theorem 5.2]{LDuke}). 
Indeed, the class of such $C^*$-algebras agrees with 
the class of all unital simple AH algebras 
with real rank zero and slow dimension growth.

\section{Central sequences}

\begin{lem}\label{appdivisible}
Let $A$ be a unital separable approximately divisible $C^*$-algebra. 
Then for any $m\in\N$, 
there exists a unital embedding of $M_m\oplus M_{m+1}$ into $A_\infty$. 
\end{lem}
\begin{proof}
Let $l=m^2-1$. 
For any finite subset $F\subset A$ and $\ep>0$, 
there exists a unital finite dimensional subalgebra $B\subset A$ 
such that every central summand of $B$ is at least $l\times l$  and 
$\lVert[a,b]\rVert<\ep$ 
for any $a\in F$ and $b\in B$ with $\lVert b\rVert\leq1$. 
It is easy to find a unital subalgebra $C$ of $B$ 
such that $C\cong M_m\oplus M_{m+1}$, 
which completes the proof. 
\end{proof}

The following is Lemma 3.6 of \cite{K98JOT}. 

\begin{lem}
Let $A$ be a unital simple AT algebra with real rank zero. 
For any finite subset $F\subset A$ and $\ep>0$, 
there exist a finite subset $G\subset A$, $\delta>0$ and $k\in\N$ 
such that the following holds. 
If $p,q\in A$ are projections satisfying $k[p]\leq[q]$, 
$\lVert[a,p]\rVert<\delta$ and $\lVert[a,q]\rVert<\delta$ 
for all $a\in G$, 
then there exists a partial isometry $v\in A$ such that 
$v^*v=p$, $vv^*\leq q$ and $\lVert[a,v]\rVert<\ep$ for all $a\in F$. 
\end{lem}

We generalize the lemma above to AH algebras. 

\begin{lem}\label{prj1}
The above lemma also holds for any unital simple AH algebra 
with slow dimension growth and real rank zero. 
\end{lem}
\begin{proof}
Let $A$ be a unital simple AH algebra 
with slow dimension growth and real rank zero and 
let $Q$ be the UHF algebra such that $K_0(Q)\cong\Q$. 
By the classification theorem (\cite{EG,D,G}), 
$A\otimes Q$ is a unital simple AT algebra with real rank zero. 
Let $F\subset A$ be a finite subset and $\ep>0$. 
Applying the lemma above to 
$\{a\otimes1\mid a\in F\}\subset A\otimes Q$ and $\ep/2>0$, 
we obtain a finite subset $G\subset A\otimes Q$, $\delta>0$ and $k\in\N$. 
We may assume $G=\{a\otimes1\mid a\in G_0\}\cup\{1\otimes b\mid b\in G_1\}$, 
where $G_0$ and $G_1$ are finite subsets of $A$ and $Q$, respectively. 
We may further assume that 
$G_0$ contains $F$ and $\delta$ is less than $\ep/2$. 
We will prove that $G_0$, $\delta$ and $4k$ meet the requirement. 

Suppose that $p,q\in A$ are non-zero projections satisfying $4k[p]\leq[q]$, 
$\lVert[a,p]\rVert<\delta$ and $\lVert[a,q]\rVert<\delta$ 
for all $a\in G_0$. 
By Lemma \ref{appdivisible}, $M_2\oplus M_3$ embeds into $A_\infty$. 
Hence there exist a projection $r$ and a partial isometry $s$ such that 
\[
r\leq q,\quad s^*s=r,\quad ss^*\leq q-r,\quad 4[r]>[q]
\]
and 
\[
\lVert[a,s]\rVert<\delta,\quad \lVert[a,r]\rVert<\delta\quad 
\forall a\in G_0. 
\]
From $4k[p]\leq[q]<4[r]$, we get $k[p]<[r]$. 
It follows that 
there exists a partial isometry $u\in A\otimes Q$ such that 
$u^*u=p\otimes1$, $uu^*\leq r\otimes1$ and 
$\lVert[a\otimes1,u]\rVert<\ep/2$ for all $a\in F$. 
We may assume that 
$u$ belongs to some $A\otimes M_m\subset A\otimes Q$. 
Put $u=(u_{i,j})_{1\leq i,j\leq m}$. 
Define $w=(w_{i,j})_{1\leq i,j\leq m+1}\in A\otimes M_{m+1}$ by 
\[
w_{i,j}=\begin{cases}
u_{i,j} & \text{if }1\leq i,j\leq m \\
su_{i,1} & \text{if }i\neq m+1\text{ and }j=m+1 \\
0 & \text{if }i=m+1. 
\end{cases}
\]
It is not so hard to see that 
$w^*w=p\otimes 1\in A\otimes M_{m+1}$ and 
$ww^*\leq q\otimes 1\in A\otimes M_{m+1}$. 
Moreover, one has $\lVert[a\otimes 1,w]\rVert<\ep$ for all $a\in F$. 
Let $v=u\oplus w\in A\otimes(M_m\oplus M_{m+1})$. 
Then $v^*v=p\otimes 1$, $vv^*\leq q\otimes 1$ and 
$\lVert[a\otimes 1,v]\rVert<\ep$ for all $a\in F$. 

By \cite{EGL} (see also \cite{D,G}), $A$ is approximately divisible. 
By Lemma \ref{appdivisible}, 
there exists a unital homomorphism from $M_m\oplus M_{m+1}$ to $A_\infty$, 
and so there exists a unital homomorphism $\pi$ 
from $A\otimes(M_m\oplus M_{m+1})$ to $A^\infty$ 
such that $\pi(a\otimes1)=a$ for $a\in A$. 
It follows that $\pi(v)^*\pi(v)=p$, $\pi(v)\pi(v)^*\leq q$ and 
$\lVert[a,\pi(v)]\rVert<\ep$ for all $a\in F$. 
\end{proof}

\begin{rem}
By using the lemma above and \cite[Theorem 4.5]{S}, 
one can show the following. 
Let $A$ be a unital simple AH algebra 
with real rank zero and slow dimension growth. 
Then for any $\alpha\in\Aut(A)$, 
there exists $\tilde\alpha\in\Aut(A)$ such that 
$\tilde\alpha$ has the Rohlin property 
in the sense of \cite[Definition 4.1]{K98JFA} and 
$\tilde\alpha$ is asymptotically unitarily equivalent to $\alpha$. 
\end{rem}

The following is a well-known fact. 
We have been unable to find a suitable reference in the literature, 
so we include a proof for completeness.

\begin{lem}\label{extrace}
Let $A$ be a unital separable $C^*$-algebra and 
let $(p_n)_n$ be a central sequence of projections. 
For any extremal trace $\tau\in T(A)$, one has 
\[
\lim_{n\to\infty}\lvert\tau(ap_n)-\tau(a)\tau(p_n)\rvert=0
\]
for all $a\in A$. 
\end{lem}
\begin{proof}
First, we deal with the case that 
there exists $\ep>0$ such that $\tau(p_n)\geq\ep$ for all $n\in\N$. 
Consider a sequence of states 
\[
\phi_n(a)=\frac{\tau(ap_n)}{\tau(p_n)}
\]
on $A$. 
Let $\psi$ be an accumulation point of $\{\phi_n\}_n$. 
Since $(p_n)_n$ is a central sequence and 
$\tau(p_n)$ is bounded from below, 
one can see that $\psi$ is a trace. 
For any positive element $a\in A$, 
it is easy to see $\phi_n(a)\leq\ep^{-1}\tau(a)$, 
and so $\psi(a)\leq\ep^{-1}\tau(a)$. 
Hence, $\psi$ is equal to $\tau$, because $\tau$ is extremal. 
We have shown that any accumulation point of $\{\phi_n\}$ is $\tau$, 
which implies $\phi_n$ converges to $\tau$. 
Therefore, $\lvert\tau(ap_n)-\tau(a)\tau(p_n)\rvert$ goes to zero. 

Next, we consider the general case. 
Fix $a\in A$. Take $\ep>0$ arbitrarily. 
We would like to show that 
$\lvert\tau(ap_n)-\tau(a)\tau(p_n)\rvert$ is less than $\ep$ 
for sufficiently large $n$. 
We may assume $\lVert a\rVert\leq1$. 
Put 
\[
C=\{n\in\N\mid\tau(p_n)\geq\ep/2\}. 
\]
If $n\notin C$, then evidently 
$\lvert\tau(ap_n)-\tau(a)\tau(p_n)\rvert$ is less than $\ep$. 
By applying the first part of the proof to $(p_n)_{n\in C}$, 
we have $\lvert\tau(ap_n)-\tau(a)\tau(p_n)\rvert<\ep$ 
for sufficiently large $n\in C$, 
thereby completing the proof. 
\end{proof}

\begin{lem}\label{prj2}
Let $A$ be a unital simple separable $C^*$-algebra 
with tracial rank zero. 
Suppose that $A$ has finitely many extremal traces. 
For any finite subset $F\subset A$ and $\ep>0$, 
there exist a finite subset $G\subset A$ and $\delta>0$ 
such that the following hold. 
If $p,q\in A$ are projections satisfying 
\[
\lVert[a,p]\rVert<\delta,\quad \lVert[a,q]\rVert<\delta\quad 
\forall a\in G
\]
and $\tau(p)+\ep<\tau(q)$ for all $\tau\in T(A)$, then 
there exists a partial isometry $v\in A$ such that 
$v^*v\leq p$, $vv^*\leq q$, 
\[
\lVert[a,v]\rVert<\ep\quad \forall a\in F
\]
and $\tau(p-v^*v)<\ep$ for all $\tau\in T(A)$. 
\end{lem}
\begin{proof}
The proof is by contradiction. 
Suppose that the assertion does not hold 
for a finite subset $F\subset A$ and $\ep>0$. 
We would have central sequences of projections $(p_n)_n$ and $(q_n)_n$ 
such that 
\[
\tau(p_n)+\ep<\tau(q_n)\quad \forall\tau\in T(A), \ n\in\N
\]
and any partial isometry does not meet the requirement. 

Use tracial rank zero to find a projection $e\in A$ and 
a finite dimensional unital subalgebra $E\subset eAe$ such that 
the following are satisfied. 
\begin{itemize}
\item For any $a\in F$, $\lVert[a,e]\rVert<\ep/4$. 
\item For any $a\in F$ there exists $b\in E$ 
such that $\lVert eae-b\rVert<\ep/4$. 
\item $\tau(1-e)<\ep$ for all $\tau\in T(A)$. 
\end{itemize}
Since $(p_n)_n$ and $(q_n)_n$ are central sequences of projections, 
we can find projections $p'_n$ and $q'_n$ in $A\cap E'$ such that 
$\lVert p_n-p'_n\rVert\to0$ and $\lVert q_n-q'_n\rVert\to0$ 
as $n\to\infty$. 
We would like to show that, 
for sufficiently large $n$, 
there exists a partial isometry $v_n\in eAe\cap E'$ such that 
$v_n^*v_n=ep'_n$ and $v_nv_n^*\leq eq'_n$. 
Let $\{e_1,e_2,\dots,e_m\}$ be 
a family of minimal central projections in $eAe\cap E'$ 
such that $e_1+e_2+\dots+e_m=e$. 
Clearly $e_i(eAe\cap E')$ is a unital simple $C^*$-algebra 
with tracial rank zero and 
the space of tracial states on $e_i(eAe\cap E')$ is naturally 
identified with $T(A)$. 
By lemma \ref{extrace}, 
for sufficiently large $n$, one has $\tau(e_ip'_n)<\tau(e_iq'_n)$ 
for all $\tau\in T(A)$ and $i=1,2,\dots,m$, 
because $A$ has finitely many extremal traces. 
Hence $[e_ip'_n]\leq[e_iq'_n]$ in $K_0(e_i(eAe\cap E'))$. 
It follows that 
there exists a partial isometry $v_n\in eAe\cap E'$ such that 
$v_n^*v_n=ep'_n$ and $v_nv_n^*\leq eq'_n$. 
Besides, $\tau(p'_n-v_n^*v_n)=\tau(p'_n(1-e))<\ep$ 
and $\lVert[a,v_n]\rVert<\ep$ for all $a\in F$. 
This is a contradiction. 
\end{proof}

By using Lemma \ref{prj2} and \ref{prj1}, 
we can show the following. 

\begin{lem}\label{prj3}
Let $A$ be a unital simple AH algebra 
with slow dimension growth and real rank zero. 
Suppose that $A$ has finitely many extremal traces. 
For any finite subset $F\subset A$ and $\ep>0$, 
there exist a finite subset $G\subset A$ and $\delta>0$ 
such that the following hold. 
If $p,q\in A$ are projections satisfying 
\[
\lVert[a,p]\rVert<\delta,\quad \lVert[a,q]\rVert<\delta\quad 
\forall a\in G
\]
and $\tau(p)+\ep<\tau(q)$ for all $\tau\in T(A)$, then 
there exists a partial isometry $u\in A$ such that 
$u^*u=p$, $uu^*\leq q$ and $\lVert[a,u]\rVert<\ep$ for all $a\in F$. 
\end{lem}
\begin{proof}
Notice that $A$ has tracial rank zero by \cite[Proposition 2.6]{LTR2}. 
Suppose that a finite subset $F\subset A$ and $\ep>0$ are given. 
By applying Lemma \ref{prj1} to $F$ and $\ep/2$, 
we obtain a finite subset $F_1\subset A$, $\ep_1>0$ and $k\in\N$. 
By applying Lemma \ref{prj2} to 
$F\cup F_1\cup F_1^*$ and $\min\{\ep_1/4,\ep/k,\ep/2\}$, 
we obtain a finite subset $G\subset A$, $\delta>0$. 
We would like to show that 
$G\cup F_1$ and $\min\{\delta,\ep_1/2\}$ meet the requirement. 
Suppose that $p,q\in A$ are projections satisfying 
\[
\lVert[a,p]\rVert<\min\{\delta,\ep_1/2\},\quad 
\lVert[a,q]\rVert<\min\{\delta,\ep_1/2\}\quad 
\forall a\in G\cup F_1
\]
and $\tau(p)+\ep<\tau(q)$ for all $\tau\in T(A)$. 
By Lemma \ref{prj2}, 
there exists a partial isometry $v\in A$ such that 
$v^*v\leq p$, $vv^*\leq q$, 
\[
\lVert[a,v]\rVert<\min\{\ep_1/4,\ep/2\}\quad 
\forall a\in F\cup F_1\cup F_1^*
\]
and $\tau(p-v^*v)<\ep/k$ for all $\tau\in T(A)$. 
Let $p'=p-v^*v$ and $q'=q-vv^*$. 
One has 
\[
\tau(q')=\tau(q-vv^*)=\tau(q)-\tau(p)+\tau(p-v^*v)>\ep, 
\]
and so $k[p']\leq[q']$. 
One also has $\lVert[a,p']\rVert<\ep_1$ and $\lVert[a,q']\rVert<\ep_1$ 
for all $a\in F_1$. 
By Lemma \ref{prj1}, we obtain a partial isometry $w\in A$ such that 
$w^*w=p'$, $ww^*\leq q'$ and $\lVert[a,w]\rVert<\ep/2$ for all $a\in F$. 
Put $u=v+w$. 
It is easy to see $u^*u=p$, $uu^*\leq q$ and 
$\lVert[a,u]\rVert<\ep$ for all $a\in F$. 
\end{proof}

Any tracial state $\tau\in T(A)$ naturally extends to 
a tracial state on $A^\omega$ and we write it by $\tau_\omega$. 

\begin{prop}\label{byproduct}
Let $A$ be a unital simple AH algebra 
with slow dimension growth and real rank zero. 
Suppose that $A$ has finitely many extremal traces. 
If $p,q\in A_\omega$ are projections satisfying 
$\tau_\omega(p)<\tau_\omega(q)$ for all $\tau\in T(A)$, 
then there exists $v\in A_\omega$ such that 
$v^*v=p$ and $vv^*\leq q$. 
In particular, $A_\omega$ satisfies 
Blackadar's second fundamental comparability question. 
\end{prop}
\begin{proof}
Let $(p_n)_n$ and $(q_n)_n$ be 
$\omega$-central sequences of projections such that 
\[
\lim_{n\to\omega}\tau(p_n)<\lim_{n\to\omega}\tau(q_n)
\]
for all $\tau\in T(A)$. 
Since $A$ has finitely many extremal traces, 
there exists $\ep>0$ such that 
\[
C=\{n\in\N\mid\tau(p_n)+\ep<\tau(q_n)\text{ for all }\tau\in T(A)\}
\in\omega. 
\]
We choose an increasing sequence 
$\{F_m\}_{m=1}^\infty$ of finite subsets of $A$ whose union is dense in $A$. 
By applying Lemma \ref{prj3} to $F_m$ and $\ep/m$, 
we obtain a finite subset $G_m\subset A$ and $\delta_m>0$. 
We may assume that 
$\{G_m\}_m$ is increasing and $\{\delta_m\}_m$ is decreasing. 
Put 
\[
C_m=\{n\in C\mid
\lVert[a,p_n]\rVert<\delta_m\text{ and }
\lVert[a,q_n]\rVert<\delta_m\text{ for all }a\in G_m\}\in\omega. 
\]
For $n\in C_m\setminus C_{m+1}$, by Lemma \ref{prj3}, 
there exists a partial isometry $u_n$ such that 
$u_n^*u_n=p_n$, $u_nu_n^*\leq q_n$ and 
$\lVert[a,u_n]\rVert<\ep/m$ for all $a\in F_m$. 
For $n\in\N\setminus C_1$, we let $u_n=0$. 
Then $(u_n)_n$ is a desired $\omega$-central sequence of partial isometries. 
\end{proof}

The following is Lemma 4.4 of \cite{K98JFA}. 

\begin{lem}
Let $A$ be a unital simple AT algebra with real rank zero. 
For any finite subset $F\subset A$ and $\ep>0$, 
there exist a finite subset $G\subset A$ and $\delta>0$ 
such that the following holds. 
If $u:[0,1]\to A$ is a path of unitaries satisfying 
$\lVert[a,u(t)]\rVert<\delta$ for all $a\in G$ and $t\in[0,1]$, 
then there exists a path of unitaries $v:[0,1]\to A$ satisfying 
\[
v(0)=u(0),\quad v(1)=u(1),\quad \lVert[a,v(t)]\rVert<\ep\quad 
\forall a\in F, \ t\in[0,1]
\]
and $\Lip(v)<5\pi+1$. 
\end{lem}

We generalize the lemma above to AH algebras. 

\begin{lem}\label{untry}
The above lemma also holds for any unital simple AH algebra 
with slow dimension growth and real rank zero, 
the Lipschitz constant being bounded by $11\pi$. 
\end{lem}
\begin{proof}
Let $A$ be a unital simple AH algebra 
with slow dimension growth and real rank zero and 
let $Q$ be the UHF algebra such that $K_0(Q)\cong\Q$. 
By the classification theorem (\cite{EG,D,G}), 
$A\otimes Q$ is a unital simple AT algebra with real rank zero. 
Let $F\subset A$ be a finite subset and $\ep>0$. 
We may assume that $F$ is contained in the unit ball of $A$. 
Applying the lemma above to 
$\{a\otimes1\mid a\in F\}\subset A\otimes Q$ and $\ep/2>0$, 
we obtain a finite subset $G\subset A\otimes Q$ and $\delta>0$. 
We may assume $G=\{a\otimes1\mid a\in G_0\}\cup\{1\otimes b\mid b\in G_1\}$, 
where $G_0$ and $G_1$ are finite subsets of $A$ and $Q$, respectively. 
We will prove that $G_0$ and $\delta$ meet the requirement. 

Suppose that $u:[0,1]\to A$ is a path of unitaries 
satisfying $\lVert[a,u(t)]\rVert<\delta$ for all $a\in G_0$ and $t\in[0,1]$. 
Choose $N\in\N$ so that 
$\lVert u(\frac{k}{N})-u(\frac{k+1}{N})\rVert<\ep/2$ 
for every $k=0,1,\dots,N{-}1$. 
Put $u_k=u(k/N)$. 
By the lemma above, we can find continuous paths $x:[0,1]\to A\otimes Q$, 
$y_k:[0,1]\to A\otimes Q$ and $z_k:[0,1]\to A\otimes Q$ 
for $k=1,2,\dots,N{-}1$ 
such that 
\[
x(0)=u_0\otimes1,\quad x(1)=u_N\otimes1,\quad 
\lVert[a\otimes1,x(t)]\rVert<\ep/2\quad \forall a\in F, \ t\in[0,1], 
\]
\[
y_k(0)=u_0\otimes1,\quad y_k(1)=u_k\otimes1,\quad 
\lVert[a\otimes1,y_k(t)]\rVert<\ep/2\quad \forall a\in F, \ t\in[0,1], 
\]
\[
z_k(0)=u_k\otimes1,\quad z_k(1)=u_N\otimes1,\quad 
\lVert[a\otimes1,z_k(t)]\rVert<\ep/2\quad \forall a\in F, \ t\in[0,1], 
\]
and $\Lip(x),\Lip(y_k),\Lip(z_k)$ are less than $5\pi+1$. 
We may assume that the ranges of $x,y_k,z_k$ are contained 
in $A\otimes M_n$ for some $M_n\subset Q$. 

Put $m=n(N-1)$. 
We would like to construct a path of unitaries 
$v:[0,1]\to A\otimes(M_m\oplus M_{m+1})$ 
such that $\Lip(v)<11\pi$, $v(0)=u_0\otimes1$, $v(1)=u_N\otimes1$ and 
$\lVert[a\otimes1,v(t)]\rVert<\ep$ for all $a\in F$ and $t\in[0,1]$. 
First, let $\tilde x:[0,1]\to A\otimes M_m$ be 
the direct sum of $N{-}1$ copies of $x:[0,1]\to A\otimes M_n$. 
Next, by using $y_1,y_2,\cdots,y_{N-1}$, 
we can find a path $\tilde y:[0,1]\to A\otimes M_{m+1}$ such that 
\[
\tilde y(0)=u_0\otimes1, 
\]
\[
\tilde y(1)=\diag(u_0,\underbrace{u_1,\cdots,u_1}_n,
\underbrace{u_2,\cdots,u_2}_n,\cdots,\underbrace{u_{N-1},\cdots,u_{N-1}}_n)
\]
\[
\lVert[a\otimes1,\tilde y(t)]\rVert<\ep/2\quad 
\forall a\in F, \ t\in[0,1]
\]
and $\Lip(\tilde y)<5\pi+1$. 
Likewise, by using $z_1,z_2,\dots,z_{N-1}$, 
we can find a path $\tilde z:[0,1]\to A\otimes M_{m+1}$ such that 
\[
\tilde z(0)=\diag(\underbrace{u_1,\cdots,u_1}_n,
\underbrace{u_2,\cdots,u_2}_n,\cdots,\underbrace{u_{N-1},\cdots,u_{N-1}}_n,u_N)
\]
\[
\tilde z(1)=u_N\otimes1, 
\]
\[
\lVert[a\otimes1,\tilde z(t)]\rVert<\ep/2\quad 
\forall a\in F, \ t\in[0,1]
\]
and $\Lip(\tilde z)<5\pi+1$. 
Since $\lVert\tilde y(1)-\tilde z(0)\rVert<\ep/2$, 
if $\ep$ is sufficiently small, 
there exists a path $w:[0,1]\to M_{m+1}$ such that 
\[
w(0)=u_0\otimes1,\quad w(1)=u_N\otimes1,\quad 
\lVert[a\otimes1,w(t)]\rVert<\ep\quad \forall a\in F, \ t\in[0,1], 
\]
and $\Lip(w)<11\pi$. 
Then $v=\tilde x\oplus w$ is the desired path. 

By \cite{EGL} (see also \cite{D,G}), $A$ is approximately divisible. 
By Lemma \ref{appdivisible}, 
there exists a unital homomorphism from $M_m\oplus M_{m+1}$ to $A_\infty$, 
and so there exists a unital homomorphism $\pi$ 
from $A\otimes(M_m\oplus M_{m+1})$ to $A^\infty$ 
such that $\pi(a\otimes1)=a$ for $a\in A$. 
It follows that 
the path $\tilde v:[0,1]\ni t\mapsto\pi(v(t))\in A^\infty$ satisfies 
\[
\tilde v(0)=u_0,\quad \tilde v(1)=u_N,\quad 
\lVert[a,\tilde v(t)]\rVert<\ep\quad \forall a\in F, \ t\in[0,1]
\]
and $\Lip(\tilde v)<11\pi$, which completes the proof. 
\end{proof}

\section{Automorphisms of AH algebras}

In this section, 
we discuss the Rohlin property of automorphisms of AH algebras. 
For $a\in A$, we define 
\[
\lVert a\rVert_2=\sup_{\tau\in T(A)}\tau(a^*a)^{1/2}. 
\]
If $A$ is simple and $T(A)$ is non-empty, 
then $\lVert\cdot\rVert_2$ is a norm. 

\begin{prop}\label{LinMatui}
Let $A$ be a unital simple separable $C^*$-algebra with tracial rank zero 
and let $\Gamma\subset\Aut(A)$ be a finite subset containing the identity. 
Suppose that there exists a sequence of projections 
$(e_n)_n$ in $A$ satisfying the following property. 
\begin{enumerate}
\item $\lVert\gamma(e_n)\gamma'(e_n)\rVert_2\to0$ 
for any $\gamma,\gamma'\in\Gamma$ such that $\gamma\neq\gamma'$. 
\item $\lVert1-\sum_{\gamma\in\Gamma}\gamma(e_n)\rVert_2\to0$. 
\item For every $a\in A$, we have 
$\lVert[a,e_n]\rVert_2\to0$. 
\end{enumerate}
Then there exists a sequence of projections 
$(f_n)_n$ in $A$ satisfying the following. 
\begin{enumerate}
\item $\lVert\gamma(f_n)\gamma'(f_n)\rVert\to0$ 
for any $\gamma,\gamma'\in\Gamma$ such that $\gamma\neq\gamma'$. 
\item $\lVert e_n-f_n\rVert_2\to0$. 
\item For every $a\in A$, we have 
$\lVert[a,f_n]\rVert\to0$. 
\end{enumerate}
\end{prop}
\begin{proof}
This is almost the same as \cite[Proposition 5.4]{LM}. 
In \cite[Proposition 5.4]{LM}, 
the finite set $\Gamma$ is assumed to be 
an orbit of a single automorphism $\gamma$ of finite order. 
The proof, however, does not need this. 
\end{proof}

The following is a variant of 
\cite[Lemma 3.1]{K95crelle} and \cite[Theorem 2.17]{OP2}. 
See \cite[Definition 2.1]{OP2} 
for the definition of the tracial Rohlin property. 

\begin{thm}\label{OsakaPhillips}
Let $A$ be a unital simple separable $C^*$-algebra with tracial rank zero. 
Suppose that $A$ has finitely many extremal traces. 
Let $\alpha$ be an automorphism of $A$ such that 
$\alpha^m$ is uniformly outer for any $m\in\N$. 
Then $\alpha$ has the tracial Rohlin property. 
\end{thm}
\begin{proof}
Let $\{\tau_1,\dots,\tau_d\}$ be the set of extremal tracial states of $A$ 
and let $(\pi_i,H_i)$ be the GNS representation associated with $\tau_i$. 
It is well-known that $\pi_i(A)''$ is a hyperfinite II$_1$-factor 
(see \cite[Lemma 2.16]{OP2}). 
Let $\rho=\bigoplus_{i=1}^d\pi_i$. 
Note that, for a bounded sequence $(a_n)_n$ in $A$, 
$\rho(a_n)$ converges to zero in the strong operator topology 
if and only if $\lVert a_n\rVert_2$ converges to zero. 
We regard $A$ as a subalgebra of 
$N=\rho(A)''\cong\bigoplus_{i=1}^d\pi_i(A)''$ and 
denote the extension of the automorphism $\alpha$ to $N$ by $\bar\alpha$. 
Let $k$ be the minimum positive integer 
such that $\tau_i\circ\alpha^k=\tau_i$ for all $i=1,2,\dots,d$. 
In the same way as \cite[Lemma 3.1]{K95crelle}, for any $l\in\N$, 
one can find a sequence $\{f^{(j)}_0,\dots,f^{(j)}_{kl-1}\}$ 
of orthogonal families of projections in $N$ such that 
$\sum_{i=0}^{kl-1}f^{(j)}_i=1$, 
\[
[a,f^{(j)}_i]\to0\quad \forall a\in A, 
\]
\[
\bar\alpha(f^{(j)}_i)-f^{(j)}_{i+1}\to0\quad \forall i=0,1,\dots,kl{-}1
\]
in the strong operator topology as $j\to\infty$, 
where $f^{(j)}_{kl}=f^{(j)}_0$. 
By \cite[Lemma 2.15]{OP2}, 
we may replace the projections $f^{(j)}_i$ with projections of $A$. 
From the proposition above, 
we can conclude that $\alpha$ has the tracial Rohlin property. 
\end{proof}

The following is a well-known fact, 
but we include the proof for the reader's convenience. 

\begin{lem}
Let $A$ be a unital separable $C^*$-algebra and 
let $\alpha\in\aInn(A)$. 
For any separable subset $C\subset A_\infty$, 
there exists a unitary $u\in A_\infty$ such that 
$uxu^*=\alpha(x)$ for all $x\in C$. 
\end{lem}
\begin{proof}
We choose an increasing sequence $\{F_n\}_{n\in\N}$ of finite subsets of $A$ 
whose union is dense in $A$. 
We can find a sequence of unitaries $(v_n)_n$ in $A$ such that 
\[
\lVert v_nav_n^*-\alpha^{-1}(a)\rVert<n^{-1}
\]
for all $a\in F_n$, because $\alpha$ is approximately inner. 
We may assume that $C$ is countable. 
Let $C=\{x_1,x_2,\dots\}$ and 
let $(x_{i,j})_j$ be a representative of $x_i$. 
There exists an increasing sequence $(m(n))_n$ of natural numbers such that 
\[
\lVert[v_n,x_{i,j}]\rVert<n^{-1}\quad \forall j\geq m(n)
\]
for any $i=1,2,\dots,n$, because $(x_{i,j})_j$ is a central sequence. 
Since $\alpha$ is in $\aInn(A)$, 
one can find a sequence of unitaries $(w_n)_n$ in $A$ such that 
\[
\lVert w_naw_n^*-\alpha(a)\rVert<n^{-1}
\]
for all $a$ in 
\[
\alpha^{-1}(F_n)\cup\{x_{i,j}\mid i=1,\dots,n, \ m(n)\leq j<m(n{+}1)\}. 
\]
For $j\in\N$, find $n\in\N$ so that $m(n)\leq j<m(n{+}1)$ and 
define a unitary $u_j$ by $u_j=w_nv_n$. 
It is easy to see 
\[
\lVert[u_j,a]\rVert<2/n\quad \forall a\in F_n
\]
and 
\[
\lVert u_jx_{i,j}u_j^*-\alpha(x_{i,j})\rVert<2/n\quad \forall i=1,\dots,n, 
\]
and so the proof is completed. 
\end{proof}

We quote the following theorem by Lin and Osaka from \cite{LO}. 
See \cite[Definition 2.4]{LO} 
for the definition of the tracial cyclic Rohlin property. 
Since we need to discuss an `equivariant version' of this theorem later, 
we would like to include the proof briefly. 

\begin{thm}[{\cite[Theorem 3.4]{LO}}]\label{LinOsaka}
Let $A$ be a unital simple separable $C^*$-algebra and 
suppose that the order on projections in $A$ is determined by traces. 
Suppose that $\alpha\in\Aut(A)$ has the tracial Rohlin property. 
If $\alpha^r$ is in $\aInn(A)$ for some $r\in\N$, 
then $\alpha$ has the tracial cyclic Rohlin property. 
\end{thm}
\begin{proof}
Take $m\in\N$ and $\ep>0$ arbitrarily. 
Let $l=r(mr+1)$. 
Since $\alpha$ has the tracial Rohlin property and 
the order on projections is determined by traces, 
there exists a central sequence of projections $(e_n)_n$ such that 
\[
\lim_{n\to\infty}\lVert e_n\alpha^i(e_n)\rVert=0\quad 
\forall i=1,2,\dots,l{-}1
\]
and 
\[
\lim_{n\to\infty}\tau(1-(e_n+\alpha(e_n)+\dots+\alpha^{l{-}1}(e_n))=0\quad 
\forall\tau\in T(A). 
\]
Let $e\in A_\infty$ be the image of $(e_n)_n$ and 
define $\tilde e$ by 
\[
\tilde e=\sum_{i=0}^{r-1}\alpha^{i(mr+1)}(e). 
\]
It follows from the lemma above that 
there exists a partial isometry $v\in A_\infty$ 
such that $v^*v=\tilde e$ and $vv^*=\alpha(\tilde e)$. 
The $C^*$-algebra $C$ generated by $v,\alpha(v),\dots,\alpha^{mr-1}(v)$ 
is isomorphic to $M_{mr+1}$ and 
its unit is $\tilde e+\alpha(\tilde e)+\dots+\alpha^{mr}(\tilde e)$. 
The rest of the proof is 
exactly the same as that of \cite[Lemma 4.3]{K95crelle} 
and we omit it. 
\end{proof}

\begin{rem}\label{Lin}
The following was shown by Lin in \cite[Theorem 3.4]{LRokhlin}. 
Let $A$ be a unital simple separable $C^*$-algebra 
with tracial rank zero. 
Suppose that $\alpha\in\Aut(A)$ has the tracial cyclic Rohlin property and 
that there exists $r\in\N$ such that 
$K_0(\alpha^r)|G=\id_G$ for some subgroup $G\subset K_0(A)$ 
for which $D_A(G)$ is dense in $D_A(K_0(A))$. 
Then $A\rtimes_\alpha\Z$ has tracial rank zero. 
\end{rem}

By using Lemma \ref{prj1}, we can show the following. 

\begin{lem}\label{ZRohlin1}
Let $A$ be a unital simple AH algebra 
with slow dimension growth and real rank zero. 
Suppose that $\alpha\in\Aut(A)$ has the tracial cyclic Rohlin property and 
that there exists $r\in\N$ such that 
$\tau\circ\alpha^r=\tau$ for any $\tau\in T(A)$. 
Then, for any $m\in\N$, there exist projections $e,f\in A_\infty$ and 
a partial isometry $v\in A_\infty$ such that 
\[
v^*v=f,\quad vv^*\leq e,\quad f+\sum_{i=0}^{mr}\alpha^i(e)=1
\]
and $\alpha^{mr+1}(e)=e$. 
\end{lem}
\begin{proof}
Suppose that we are given $m\in\N$. 
Let $l=r(mr+1)$. 
Since $\alpha$ has the tracial cyclic Rohlin property, 
we can find central sequences of projections $(e_n)_n$ and $(f_n)_n$ 
such that
\[
f_n+\sum_{i=0}^{l-1}\alpha^i(e_n)\to1,\quad 
e_n-\alpha^l(e_n)\to0,\quad \sup_{\tau\in T(A)}\tau(f_n)\to0
\]
as $n\to\infty$. 
There exists a central sequence of projections $(\tilde{e}_n)_n$ such that 
\[
\lim_{n\to\infty}\tilde{e}_n-\sum_{i=0}^{r-1}\alpha^{i(mr+1)}(e_n)=0. 
\]
Then 
\[
f_n+\sum_{i=0}^{mr}\alpha^i(\tilde{e}_n)\to1,\quad 
\tilde{e}_n-\alpha^{mr+1}(\tilde{e}_n)\to0
\]
as $n\to\infty$. 
It is also easy to see 
$\tau(\alpha(\tilde{e}_n))=\tau(\tilde{e}_n)$ for all $\tau\in T(A)$, 
and so $\tau(\tilde{e}_n)$ goes to $(mr+1)^{-1}$ for all $\tau\in T(A)$. 
Therefore, by Lemma \ref{prj1}, 
one can find a central sequence of partial isometries $(v_n)_n$ such that 
$v_n^*v_n=f_n$ and $v_nv_n^*\leq\tilde{e}_n$ for sufficiently large $n$, 
which completes the proof. 
\end{proof}

By using the lemma above, we can show the following theorem. 

\begin{thm}\label{ZRohlin2}
Let $A$ be a unital simple AH algebra 
with slow dimension growth and real rank zero. 
Suppose that $\alpha\in\Aut(A)$ has the tracial cyclic Rohlin property. 
If there exists $r\in\N$ such that 
$\tau\circ\alpha^r=\tau$ for any $\tau\in T(A)$, 
then $\alpha$ has the Rohlin property. 
\end{thm}
\begin{proof}
Suppose that we are given $M\in\N$. 
Choose a natural number $m\in\N$ 
so that $m\geq M$ and $m\equiv1\pmod{r}$. 
Let $k,l$ be sufficiently large natural numbers 
satisfying $k\equiv l\equiv1\pmod{r}$. 
By the lemma above, 
we can find projections $e,f\in A_\infty$ and 
a partial isometry $v\in A_\infty$ such that 
\[
v^*v=f,\quad vv^*\leq e,\quad f+\sum_{i=0}^{klm-1}\alpha^i(e)=1
\]
and $\alpha^{klm}(e)=e$. 
Define $\tilde{e},w\in A_\infty$ by 
\[
\tilde{e}=\sum_{i=0}^{k-1}\alpha^{ilm}(e)\quad\text{ and }\quad 
w=\frac{1}{\sqrt{k}}\sum_{i=0}^{k-1}\alpha^{ilm}(v). 
\]
Then $\tilde{e}$ is a projection and $w$ is a partial isometry satisfying 
\[
f+\sum_{i=0}^{lm-1}\alpha^i(\tilde{e})=1,\quad 
\alpha^{lm}(\tilde{e})=\tilde{e}
\]
and 
\[
w^*w=f,\quad ww^*\leq \tilde{e},\quad 
\lVert\alpha^{lm}(w)-w\rVert\leq\frac{2}{\sqrt{k}}. 
\]
Let $D$ be the $C^*$-algebra generated 
by $w,\alpha(w),\dots,\alpha^{lm-1}(w)$. 
Then $D$ is isomorphic to $M_{lm+1}$ and 
the unit $1_D$ of $D$ is equal to $f+ww^*+\dots+\alpha^{lm-1}(ww^*)$. 
From the spectral property of $\alpha$ restricted to $D$, 
if $k$ and $l$ are sufficiently large, 
we can obtain projections $p_0,\dots,p_{m-1},q_0,\dots,q_m$ of $D$ such that 
\[
\sum_{i=1}^{m-1}p_i+\sum_{i=1}^mq_i=1_D,\quad 
\alpha(p_i)\approx p_{i+1},\quad \alpha(q_i)\approx q_{i+1}, 
\]
where $p_m=p_0$ and $q_{m+1}=q_0$. 
We define projections $p'_i$ in $A_\infty$ by 
\[
p'_i=p_i+\sum_{j=0}^{l-1}\alpha^{i+jm}(\tilde{e}-ww^*). 
\]
Then the projections $p'_0,\dots,p'_{m-1},q_0,\dots,q_m$ meet the requirement. 
See \cite{K95crelle,K96JFA} for details. 
\end{proof}

Combining the theorems above, we obtain the following theorem 
which is a generalization of \cite[Theorem 2.1]{K98JOT}. 

\begin{thm}\label{ZRohlintype}
Let $A$ be a unital simple AH algebra 
with slow dimension growth and real rank zero 
and let $\alpha\in\Aut(A)$. 
Suppose that $A$ has finitely many extremal traces and 
that $\alpha^r$ is approximately inner for some $r\in\N$. 
Then the following are equivalent. 
\begin{enumerate}
\item $\alpha$ has the Rohlin property. 
\item $\alpha^m$ is uniformly outer for any $m\in\N$. 
\end{enumerate}
\end{thm}

We can also generalize \cite[Theorem 5.1]{K98JFA} 
by using Lemma \ref{untry} instead of \cite[Lemma 4.4]{K98JFA}. 

\begin{thm}\label{Zcc}
Let $A$ be a unital simple AH algebra 
with slow dimension growth and real rank zero. 
If $\alpha,\beta\in\Aut(A)$ have the Rohlin property and 
$\alpha$ is asymptotically unitarily equivalent to $\beta$, then 
there exist $\mu\in\aInn(A)$ and a unitary $u\in A$ such that 
\[
\Ad u\circ\alpha=\mu\circ\beta\circ\mu^{-1}. 
\]
\end{thm}

The proof is similar to that of \cite[Theorem 5.1]{K98JFA} 
and we omit it. 

As an application of the theorems above, 
we can show the following, which will be used in Section 6. 

\begin{lem}\label{equivhomotopy}
Let $A$ be a unital simple AF algebra with finitely many extremal traces 
and let $\alpha$ be an approximately inner automorphism of $A$ such that 
$\alpha^m$ is uniformly outer for all $m\in\N$. 
For any finite subset $F\subset A$ and $\ep>0$, 
there exist a finite subset $G\subset A$ and $\delta>0$ 
satisfying the following. 
If $u:[0,1]\to A$ is a path of unitaries such that 
\[
\lVert[a,u(t)]\rVert<\delta\quad\text{and}\quad 
\lVert u(t)-\alpha(u(t))\rVert<\delta\quad 
\forall a\in G,\ t\in[0,1], 
\]
then there exists a path of unitaries $v:[0,1]\to A$ such that 
\[
v(0)=u(0),\quad v(1)=u(1),\quad 
\lVert[a,v(t)]\rVert<\ep,\quad 
\lVert v(t)-\alpha(v(t))\rVert<\ep,\quad 
\forall a\in F,\ t\in[0,1]
\]
and $\Lip(v)<2\pi$. 
\end{lem}
\begin{proof}
Let $x_n$ be a unitary of $M_n(\C)$ such that 
$\Sp(x_n)=\{\omega^k\mid k=0,1,\dots,n{-}1\}$, 
where $\omega=\exp(2\pi\sqrt{-1}/n)$. 
One can find an increasing sequence $\{A_n\}_{n=1}^\infty$ of 
unital finite dimensional subalgebras of $A$ such that 
$\bigcup_n A_n$ is dense in $A$ and 
there exists a unital embedding $\pi_n:M_n\oplus M_{n+1}\to A_{n+1}\cap A_n'$. 
Let $y_n=\pi_n(x_n\oplus x_{n+1})$. 
Define an automorphism $\sigma$ of $A$ 
by $\sigma=\lim_{n\to\infty}\Ad(y_1y_2\dots y_n)$. 
Then $\sigma$ is approximately inner and 
$\sigma^m$ is uniformly outer for all $m\in\N$. 

We would like to show that the assertion holds for $\sigma$. 
Suppose that we are given $F\subset A$ and $\ep>0$. 
Without loss of generality, 
we may assume that there exists $n\in\N$ such that 
$F$ is contained in the unit ball of $A_n$. 
Applying \cite[Lemma 4.2]{KM} to $\ep/2$, 
we obtain a positive real number $\delta_1>0$. 
We may assume $\delta_1$ is less than $\min\{2,\ep\}$. 
Choose a finite subset $G\subset A$ and $\delta_2>0$ so that 
if $z:[0,1]\to A$ is a path of unitaries such that 
$\lVert[a,z(t)]\rVert<\delta_2$ for all $a\in G$ and $t\in[0,1]$, 
then there exists a path of unitaries $\tilde z:[0,1]\to A\cap A_n'$ 
such that $\lVert z(t)-\tilde z(t)\rVert<\delta_1/6$. 
Let $\delta=\min\{\delta_1/6,\delta_2\}$. 
Suppose that $u:[0,1]\to A$ is a path of unitaries such that 
\[
\lVert[a,u(t)]\rVert<\delta\quad\text{and}\quad 
\lVert u(t)-\sigma(u(t))\rVert<\delta\quad 
\forall a\in G,\ t\in[0,1]. 
\]
By the choice of $\delta$, 
we can find $\tilde u:[0,1]\to A\cap A_n'$ 
such that $\lVert u(t)-\tilde u(t)\rVert<\delta_1/6$. 
We may assume that there exists $m>n$ such that 
the range of $\tilde u$ is contained in $A_m$. 
Put $y=y_ny_{n+1}\dots y_{m-1}\in A_m\cap A_n'$. 
Then 
\[
\lVert[y,\tilde u(t)]\rVert=\lVert\tilde u(t)-\sigma(\tilde u(t))\rVert
<\lVert u(t)-\sigma(u(t))\rVert+\delta_1/3<\delta+\delta_1/3\leq\delta_1/2
\]
for every $t\in[0,1]$. 
Hence $\lVert[y,\tilde u(t)\tilde u(0)^*]\rVert$ is less than $\delta_1$. 
It follows from \cite[Lemma 4.2]{KM} that 
one can find a path of unitaries $w:[0,1]\to A_m\cap A_n'$ such that 
\[
w(0)=1,\quad w(1)=\tilde u(1)\tilde u(0)^*,\quad \Lip(w)\leq\pi+\ep
\]
and $\lVert[y,w(t)]\rVert<\ep/2$ for every $t\in[0,1]$. 
Note that $yw(t)y^*$ is equal to $\sigma(w(t))$. 
By perturbing $w(t)u(0)$ a little bit, 
the required $v:[0,1]\to A$ is obtained. 

Suppose that 
$\alpha$ is an approximately inner automorphism of $A$ such that 
$\alpha^m$ is uniformly outer for all $m\in\N$. 
By Theorem \ref{ZRohlintype} and Theorem \ref{Zcc}, 
there exist $\mu\in\Aut(A)$ and a unitary $u\in A$ such that 
$\Ad u\circ\alpha=\mu\circ\sigma\circ\mu^{-1}$. 
Moreover, one can choose $u$ arbitrarily close to $1$, 
because $A$ is AF (see \cite{HO,K98JFA}). 
Therefore the assertion also holds for $\alpha$. 
\end{proof}

\begin{rem}\label{modelaction}
In the proof of the lemma above, it is easily seen that 
(the $\Z$-action generated by) $\sigma$ is asymptotically representable 
(\cite[Definition 2.2]{IM}). 
Hence the automorphism $\alpha$ stated in the lemma above is 
also asymptotically representable. 
Besides, it is not so hard to see that 
the crossed product $C^*$-algebra $A\rtimes_\sigma\Z$ is 
a unital simple AT algebra with real rank zero. 
Therefore, $A\rtimes_\alpha\Z$ is 
a unital simple AT algebra with real rank zero, too. 
\end{rem}

\section{The Rohlin property of $\Z^2$-actions on AF algebras}

In this section, we would like to show that 
certain $\Z^2$-actions on an AF algebra have the Rohlin property. 
This is a generalization of Nakamura's theorem \cite[Theorem 3]{N1}. 

Throughout this section, we keep the following setting. 
Let $A$ be a unital simple separable $C^*$-algebra with tracial rank zero 
and suppose that $A$ has a unique tracial state $\tau$. 
Suppose that 
automorphisms $\alpha,\beta\in\Aut(A)$ and a unitary $w\in A$ satisfy 
\[
\beta\circ\alpha=\Ad w\circ\alpha\circ\beta
\]
and 
$\alpha^m\circ\beta^n$ is uniformly outer 
for all $(m,n)\in\Z^2\setminus\{(0,0)\}$. 
We remark that $\alpha$ and $\beta$ induce a $\Z^2$-action on $A_\infty$. 

\begin{lem}\label{Z2Rohlin1}
For any $m_1,m_2\in\N$, 
there exists a central sequence of projections $(e_n)_n$ in $A$ such that 
\[
\lim_{n\to\infty}\tau(e_n)=\frac{1}{m_1m_2}
\]
and 
\[
\lim_{n\to\infty}\lVert\beta^j(\alpha^i(e_n))\beta^l(\alpha^k(e_n))\rVert=0
\]
for all $(i,j)\neq(k,l)$ 
in $\{(i,j)\mid0\leq i\leq m_1{-}1, \ 0\leq j\leq m_2{-}1\}$. 
\end{lem}
\begin{proof}
Set $I=\{(i,j)\mid0\leq i\leq m_1{-}1, \ 0\leq j\leq m_2{-}1\}$. 
Let $(\pi_\tau,H_\tau)$ be the GNS representation associated with $\tau$. 
It is well-known that $\pi_\tau(A)''$ is a hyperfinite II$_1$-factor 
(see \cite[Lemma 2.16]{OP2}). 
For $(i,j)\in\Z^2$, we put $\phi_{(i,j)}=\alpha^i\circ\beta^j$. 
Then $\phi:\Z^2\to\Aut(A)$ is a cocycle action of $\Z^2$. 
We denote its extension to $\pi_\tau(A)''$ by $\bar\phi$. 
Since $\alpha^m\circ\beta^n$ is uniformly outer 
for all $(m,n)\in\Z^2\setminus\{(0,0)\}$, 
$\bar\phi$ is an outer cocycle action of $\Z^2$ on $\pi_\tau(A)''$. 
It follows from \cite{O} that 
there exists a sequence of projections $(e_n)_n$ in $\pi_\tau(A)''$ such that 
\[
\sum_{(i,j)\in I}\beta^j(\alpha^i(e_n))\to1,\quad 
[x,e_n]\to0\quad \forall x\in\pi_\tau(A)''
\]
and 
\[
\beta^j(\alpha^i(e_n))\beta^l(\alpha^k(e_n))\to0\quad 
\forall (i,j),(k,l)\in I\text{ with }(i,j)\neq(k,l)
\]
in the strong operator topology as $n\to\infty$. 
By \cite[Lemma 2.15]{OP2}, 
we may replace $e_n$ with projections in $A$. 
Applying Proposition \ref{LinMatui} to 
$\Gamma=\{\beta^j\circ\alpha^i\mid(i,j)\in I\}$, we obtain the conclusion. 
\end{proof}

Let $(e_n)_n$ be the projections as in the lemma above. 
If $\alpha^r$ is in $\aInn(A)$ for some $r\in\N$ and 
$m_1$ is large enough, 
then we can construct a central sequence of projections $(e'_n)_n$ in $A$ 
such that 
\[
\lim_{n\to\infty}\lVert e'_n(e_n+\alpha(e_n)+\dots+\alpha^{m_1-1}(e_n))
-e'_n\rVert=0, 
\]
\[
\alpha(e'_n)\approx e'_n\quad\text{and}\quad 
\tau(e'_n)\approx m_1\tau(e_n), 
\]
by using the arguments in \cite[Lemma 3.1]{K95crelle} 
(see also \cite[Lemma 6]{N1}). 
Consequently, we get the following. 

\begin{lem}\label{Z2Rohlin2}
If $\alpha^r$ is in $\aInn(A)$ for some $r\in\N$, 
then for any $m\in\N$, 
there exists a central sequence of projections $(e_n)_n$ in $A$ such that 
\[
\lim_{n\to\infty}\tau(e_n)=\frac{1}{m},\quad 
\lim_{n\to\infty}\lVert e_n-\alpha(e_n)\rVert=0
\]
and 
\[
\lim_{n\to\infty}\lVert e_n\beta^j(e_n)\rVert=0
\]
for all $j=1,2,\dots,m{-}1$. 
\end{lem}

Our next task is 
to achieve the cyclicity condition $\beta^m(e_n)\approx e_n$. 

\begin{lem}\label{Z2Rohlin3}
Suppose that $A$ is AF. 
Suppose that 
a projection $e\in A_\infty$ and a partial isometry $u\in A_\infty$ satisfy 
$e=\alpha(e)$ and $e=u^*u=uu^*$. 
Then there exists a partial isometry $w\in A_\infty$ such that 
$w^*w=ww^*=e$ and $u=w^*\alpha(w)$. 
\end{lem}
\begin{proof}
Theorem \ref{ZRohlintype} tells us that 
$\alpha$ possesses the Rohlin property. 
We can modify the standard argument 
deducing stability from the Rohlin property (see \cite{HO,EK}) 
and apply it to the unitary $u+(1-e)$. 
We leave the details to the readers. 
\end{proof}

\begin{lem}\label{Z2Rohlin4}
Suppose that either of the following holds. 
\begin{enumerate}
\item $A$ is AF, $\alpha^r$ is approximately inner for some $r\in\N$ 
and $\beta^s$ is approximately inner for some $s\in\N$. 
\item $\alpha^r$ is approximately inner for some $r\in\N$ and 
there exist a natural number $s\in\N$ and 
a sequence of unitaries $(u_n)_n$ in $A$ such that 
\[
\lim_{n\to\infty}\lVert u_n-\alpha(u_n)\rVert=0\quad\text{ and }\quad 
\lim_{n\to\infty}\lVert u_nau_n^*-\beta^s(a)\rVert=0\quad \forall a\in A. 
\]
\end{enumerate}
Then for any $m\in\N$, 
there exists a central sequence of projections $(e_n)_n$ such that 
\[
\lim_{n\to\infty}\tau(e_n)=\frac{1}{m},\quad 
\lim_{n\to\infty}\lVert e_n-\alpha(e_n)\rVert=0,\quad 
\lim_{n\to\infty}\lVert e_n\beta^j(e_n)\rVert=0
\]
for all $j=1,2,\dots,m{-}1$ and 
\[
\lim_{n\to\infty}\lVert e_n-\beta^m(e_n)\rVert=0. 
\]
\end{lem}
\begin{proof}
Choose a large natural number $l$ such that $l\equiv1\pmod{s}$. 
By using Lemma \ref{Z2Rohlin2} and 
the assumption that $\beta^s$ is in $\aInn(A)$ for some $s\in\N$, 
one can find a projection $e\in A_\infty$ and 
a partial isometry $v\in A_\infty$ such that 
\[
e=\alpha(e),\quad v^*v=e,\quad vv^*=\beta(e)\quad\text{ and }\quad
e\beta^j(e)=0\quad \forall j=1,2,\dots,l{-}1
\]
in the same way as the proof of Theorem \ref{LinOsaka}. 
Moreover, we have 
\[
\lim_{n\to\infty}\tau(e_n)=l^{-1}, 
\]
where $(e_n)_n$ is a representative sequence of $e$ consisting of projections. 
Note that $\beta^j(e)$ is fixed by $\alpha$, 
because $e$ is a central sequence. 
In the case (2), clearly we may further assume $v=\alpha(v)$. 
In the case (1), 
the lemma above applies to $v^*\alpha(v)$ and yields $w\in A_\infty$ 
satisfying $w^*w=e$, $ww^*=e$ and $v^*\alpha(v)=w^*\alpha(w)$. 
By replacing $v$ with $vw^*$, we get $v=\alpha(v)$, too. 
Then the conclusion follows from 
exactly the same argument as \cite[Lemma 4.3]{K95crelle}. 
\end{proof}

\begin{thm}\label{Z2Rohlin5}
Suppose that the conclusion of Lemma \ref{Z2Rohlin4} holds. 
Then for any $m\in\N$, 
there exist projections $e$ and $f$ in $A_\infty$ such that 
\[
\alpha(e)=e,\quad \alpha(f)=f,\quad \beta^m(e)=e,\quad \beta^{m+1}(f)=f
\]
and 
\[
\sum_{i=0}^{m-1}\beta^i(e)+\sum_{j=0}^m\beta^j(f)=1. 
\]
\end{thm}
\begin{proof}
Let $(e_n)_n$ be the central sequence of projections 
obtained in Lemma \ref{Z2Rohlin4}. 
Define 
\[
f_n=1-\sum_{j=0}^{m-1}\beta^j(e_n). 
\]
There exists a sequence of unitaries $(u_n)_n$ in $A$ such that 
$u_n\to1$ as $n\to\infty$ and 
$u_n\alpha(e_n)u_n^*=e_n$ for sufficiently large $n$. 
The $\Z$-action on $e_nAe_n$ generated by $\Ad u_n\circ\alpha$ is 
uniformly outer, 
and so it has the tracial Rohlin property by Theorem \ref{OsakaPhillips} 
(or \cite[Theorem 2.17]{OP2}). 
It follows that, for any $k\in\N$, 
there exists a central sequence of projections $(\tilde e_n)_n$ such that 
\[
\tilde e_n\leq e_n,\quad 
\lim_{n\to\infty}\tau(\tilde e_n)=1/mk,\quad\text{and}\quad 
\lim_{n\to\infty}\lVert\tilde e_n\alpha^i(\tilde e_n)\rVert=0\quad 
\forall i=1,2,\dots,k{-}1. 
\]
Let $e,f,\tilde e\in A_\infty$ be the images of 
$(e_n)_n,(f_n)_n,(\tilde e_n)_n$, respectively. 
By Lemma \ref{prj1}, 
there exists a partial isometry $v$ such that 
$v^*v=f$ and $vv^*\leq \tilde e$. 
We define a partial isometry $\tilde v\in A_\infty$ by 
\[
\tilde v=\frac{1}{\sqrt{k}}\sum_{i=0}^{k-1}\alpha^i(v). 
\]
Then one has 
\[
\tilde v^*\tilde v=f,\quad \tilde v\tilde v^*\leq e\quad\text{and}\quad 
\lVert\tilde v-\alpha(\tilde v)\rVert<2/\sqrt{k}. 
\]
By a standard trick on central sequences, 
we may assume $\alpha(\tilde v)=\tilde v$. 
Thus, we have obtained 
the $\alpha$-invariant version of the conclusion of Lemma \ref{ZRohlin1}. 
We can complete the proof by the same argument as in Theorem \ref{ZRohlin2}. 
\end{proof}

The following is a generalization of \cite[Theorem 3]{N1}. 

\begin{cor}\label{Z2Rohlin6}
Let $\phi:\Z^2\curvearrowright A$ be a $\Z^2$-action 
on a unital simple AF algebra $A$ with unique trace. 
When $\phi_{(r,0)}$ and $\phi_{(0,s)}$ are approximately inner 
for some $r,s\in\N$, 
the following are equivalent. 
\begin{enumerate}
\item $\phi$ has the Rohlin property. 
\item $\phi$ is uniformly outer. 
\end{enumerate}
\end{cor}
\begin{proof}
This immediately follows 
from Theorem \ref{Z2Rohlin5} and \cite[Remark 2]{N1} 
(see also \cite[Remark 2.2]{M08}). 
\end{proof}

The next corollary also follows from Theorem \ref{Z2Rohlin5} immediately, 
because condition (2) of Lemma \ref{Z2Rohlin4} is satisfied in this case. 
See \cite[Definition 2.2]{IM} 
for the definition of approximate representability. 

\begin{cor}
Let $\phi:\Z^2\curvearrowright A$ be 
an approximately representable $\Z^2$-action 
on a unital simple AH algebra $A$ 
with real rank zero and slow dimension growth. 
Suppose that $A$ has a unique trace. 
Then the following are equivalent. 
\begin{enumerate}
\item $\phi$ has the Rohlin property. 
\item $\phi$ is uniformly outer. 
\end{enumerate}
\end{cor}

\section{Classification of certain $\Z^2$-actions on AF algebras}

In this section, we will show a classification result 
of a certain class of $\Z^2$-actions on unital simple AF algebras. 
We freely use the terminology and notation 
introduced in \cite[Definition 2.1]{IM}. 
For an automorphism $\alpha$ of a $C^*$-algebra $A$, 
we write the crossed product $C^*$-algebra $A\rtimes_\alpha\Z$ 
by $C^*(A,\alpha)$ and the implementing unitary by $\lambda_\alpha$. 
The mapping torus $M(A,\alpha)$ is defined by 
\[
M(A,\alpha)=\{f\in C([0,1],A)\mid\alpha(f(0))=f(1)\}. 
\]
When $A$ is an AF algebra, 
`$KK$-triviality' of $\alpha\in\Aut(A)$ is equivalent to $K_0(\alpha)=\id$, 
and also equivalent to $\alpha$ being approximately inner. 

The following theorem is a $\Z$-equivariant version of Theorem \ref{Zcc}. 
Let $A$ be a unital simple AF algebra with unique trace and 
let $\alpha\in\aInn(A)$. 
Let $\Aut_\T(C^*(A,\alpha))$ denote 
the set of all automorphisms of $C^*(A,\alpha)$ 
commuting with the dual action $\hat\alpha$. 
For $i=1,2$, we suppose that 
an automorphism $\beta_i\in\Aut(A)$ and 
a unitary $w_i\in A$ are given and satisfy 
\[
\beta_i\circ\alpha=\Ad w_i\circ\alpha\circ\beta_i. 
\]
Then $\beta_i$ extends to $\tilde\beta_i\in\Aut_\T(C^*(A,\alpha))$ 
by setting $\tilde\beta_i(\lambda_\alpha)=w_i\lambda_\alpha$. 
Suppose further that $\alpha^m\circ\beta_i^n$ is uniformly outer 
for all $(m,n)\in\Z^2\setminus\{(0,0)\}$ 
and that $\beta_i^{s_i}$ is approximately inner for some $s_i\in\N$. 

\begin{thm}\label{equivEK}
In the setting above, if $\tilde\beta_1$ and $\tilde\beta_2$ are 
asymptotically unitarily equivalent, then 
there exist an approximately inner automorphism $\mu\in\Aut_\T(C^*(A,\alpha))$ 
and a unitary $v\in A$ such that 
$\mu|A$ is also approximately inner and 
\[
\mu\circ\tilde\beta_1\circ\mu^{-1}=\Ad v\circ\tilde\beta_2. 
\]
\end{thm}
\begin{proof}
We can apply the argument of \cite[Theorem 5]{N2} 
to $\tilde\beta_1$ and $\tilde\beta_2$ 
in a similar fashion to \cite[Theorem 4.11]{IM}. 
By Remark \ref{modelaction}, 
(the $\Z$-action generated by) $\alpha$ is asymptotically representable. 
Then \cite[Theorem 4.8]{IM} implies that 
$\tilde\beta_1$ and $\tilde\beta_2$ are 
$\T$-asymptotically unitarily equivalent. 
Moreover, by Theorem \ref{Z2Rohlin5}, 
we can find Rohlin projections for $\tilde\beta_i$ 
in the fixed point algebra $(A_\infty)^\alpha$. 
Hence, by using Lemma \ref{equivhomotopy} instead of \cite[Theorem 7]{N2},
the usual intertwining argument shows the statement.
\end{proof}

Let us recall the $\OrderExt$ invariant introduced in \cite{KK}. 
Let $G_0,G_1,F$ be abelian groups and let $D:G_0\to F$ be a homomorphism. 
When 
\[
\begin{CD}
\xi:\quad 0 @>>> G_0 @>\iota>> E_\xi @>q>> G_1 @>>> 0
\end{CD}
\]
is exact, $R$ is in $\Hom(E_\xi,F)$ and $R\circ\iota=D$, 
the pair $(\xi,R)$ is called an order-extension. 
Two order-extensions $(\xi,R)$ and $(\xi',R')$ are equivalent 
if there exists an isomorphism $\theta:E_\xi\to E_{\xi'}$ such that 
$R=R'\circ\theta$ and 
\[
\begin{CD}
\xi:\quad 0 @>>> G_0 @>>> E_\xi @>>> G_1 @>>> 0 \\
@. @| @VV\theta V @| @. \\
\xi':\quad 0 @>>> G_0 @>>> E_{\xi'} @>>> G_1 @>>> 0
\end{CD}
\]
is commutative. 
Then $\OrderExt(G_1,G_0,D)$ consists of 
equivalence classes of all order-extensions. 
As shown in \cite{KK}, 
$\OrderExt(G_1,G_0,D)$ is equipped with an abelian group structure. 
The map sending $(\xi,R)$ to $\xi$ induces 
a homomorphism from $\OrderExt(G_1,G_0,D)$ onto $\Ext(G_1,G_0)$. 

Let $B$ be a unital $C^*$-algebra with $T(B)$ non-empty. 
We denote by $\Aut_0(B)$ the set of all automorphisms $\gamma$ of $B$ 
such that $K_0(\gamma)=K_1(\gamma)=\id$ and 
$\tau\circ\gamma=\tau$ for all $\tau\in T(B)$. 
When $B$ is a unital simple AT algebra with real rank zero, 
$\Aut_0(B)$ equals $\aInn(B)$. 
Let $D_B:K_0(B)\to\Aff(T(B))$ denote the dimension map 
defined by $D_B([p])(\tau)=\tau(p)$. 
As described in \cite{KK}, 
there exist natural homomorphisms 
\[
\tilde\eta_0:\Aut_0(B)\to\OrderExt(K_1(B),K_0(B),D_B)
\]
and 
\[
\eta_1:\Aut_0(B)\to\Ext(K_0(B),K_1(B)). 
\]
The following is the main result of \cite{KK}. 
See \cite{M02,Lasymp} for further developments. 

\begin{thm}[{\cite[Theorem 4.4]{KK}}]\label{KishiKum}
Suppose that $B$ is a unital simple AT algebra with real rank zero. 
Then the homomorphism 
\[
\tilde\eta_0\oplus\eta_1:\aInn(B)\to
\OrderExt(K_1(B),K_0(B),D_B)\oplus\Ext(K_0(B),K_1(B))
\]
is surjective and 
its kernel equals the set of all asymptotically inner automorphisms of $B$. 
\end{thm}

By using this $\OrderExt$ invariant, 
we introduce an invariant of certain $\Z^2$-actions as follows. 
Let $A$ be a unital simple AF algebra and 
let $\phi:\Z^2\curvearrowright A$ be an action of $\Z^2$ on $A$. 
Suppose that $\phi$ is uniformly outer and locally $KK$-trivial 
(i.e. locally approximately inner). 
We write $B=C^*(A,\phi_{(1,0)})$. 
Then $\phi_{(0,1)}$ extends to $\tilde\phi_{(0,1)}\in\Aut(B)$ 
by setting 
$\tilde\phi_{(0,1)}(\lambda_{\phi_{(1,0)}})=\lambda_{\phi_{(1,0)}}$. 
Let $\iota:A\to B=C^*(A,\phi_{(1,0)})$ be the canonical inclusion. 
One can check the following immediately. 
\begin{itemize}
\item $K_0(\iota)$ is an isomorphism from $K_0(A)$ to $K_0(B)$. 
\item The connecting map $\partial:K_1(B)\to K_0(A)$ 
in the Pimsner-Voiculescu exact sequence is an isomorphism and 
$\partial^{-1}([p])=[\lambda_{\phi_{(1,0)}}\iota(p)+\iota(v(1-p))]$ 
for any projection $p\in A$, 
where $v$ is a unitary of $A$ satisfying $vpv^*=\phi_{(1,0)}(p)$. 
\item The map $\iota^*:T(B)\to T(A)$ sending $\tau$ to $\tau\circ\iota$ 
is an isomorphism and 
satisfies $D_B(K_0(\iota)(x))(\tau)=D_A(x)(\iota^*(\tau))$ 
for $x\in K_0(A)$ and $\tau\in T(B)$. 
\end{itemize}
From these properties, we can obtain a natural isomorphism 
\[
\zeta_{\phi_{(1,0)}}:\OrderExt(K_1(B),K_0(B),D_B)\to
\OrderExt(K_0(A),K_0(A),D_A). 
\]
In addition, it is easy to see 
$K_0(\tilde\phi_{(0,1)})=K_1(\tilde\phi_{(0,1)})=\id$ and 
$\tau\circ\tilde\phi_{(0,1)}=\tau$ for all $\tau\in T(B)$, 
that is, $\tilde\phi_{(0,1)}$ belongs to $\Aut_0(B)$. 

\begin{lem}\label{eta1=0}
In the setting above, 
$\eta_1(\tilde\phi_{(0,1)})\in\Ext(K_0(B),K_1(B))$ is zero. 
\end{lem}
\begin{proof}
There exists a natural commutative diagram 
\[
\begin{CD}
0 @>>> C_0((0,1),B) @>>> M(B,\tilde\phi_{(0,1)}) @>>> B @>>> 0 \\
@. @AAA @AAA @AA\iota A \\
0 @>>> C_0((0,1),A) @>>> M(A,\phi_{(0,1)}) @>>> A @>>> 0, 
\end{CD}
\]
where the horizontal sequences are exact. 
From the naturality of the six-term exact sequence, 
we obtain the commutative diagram 
\[
\begin{CD}
0 @>>> K_1(B) @>>> K_0(M(B,\tilde\phi_{(0,1)})) @>>> K_0(B) @>>> 0 \\
@. @AAA @AAA @AA K_0(\iota)A \\
0 @>>> K_1(A) @>>> K_0(M(A,\phi_{(0,1)})) @>>> K_0(A) @>>> 0, 
\end{CD}
\]
where the horizontal sequences are exact. 
Since $K_1(A)$ is zero and $K_0(\iota)$ is an isomorphism, 
we can conclude $\eta_1(\tilde\phi_{(0,1)})=0$. 
\end{proof}

\begin{df}
In the setting above, 
we define our invariant $[\phi]$ by 
\[
[\phi]=\zeta_{\phi_{(1,0)}}(\tilde\eta_0(\tilde\phi_{(0,1)}))
\in\OrderExt(K_0(A),K_0(A),D_A). 
\]
\end{df}

\begin{prop}\label{invariance}
Let $\phi,\psi:\Z^2\curvearrowright A$ be 
uniformly outer, locally $KK$-trivial $\Z^2$-actions 
on a unital simple AF algebra $A$. 
If $\phi$ and $\psi$ are $KK$-trivially cocycle conjugate, 
then $[\phi]=[\psi]$. 
\end{prop}
\begin{proof}
For $\mu\in\aInn(A)$, 
it is straightforward to see that 
the $\Z^2$-action $\mu\circ\phi\circ\mu^{-1}$ has 
the same invariant as $\phi$. 
Hence, it suffices to show $[\phi]=[\phi^u]$ 
for any $\phi$-cocycle $\{u_n\}_{n\in\Z^2}$. 
Define an isomorphism $\pi$ 
from $C^*(A,\phi_{(1,0)})$ to $C^*(A,\phi^u_{(1,0)})$ by 
\[
\pi(\lambda_{\phi_{(1,0)}})=u_{(1,0)}^*\lambda_{\phi^u_{(1,0)}}
\quad\text{and}\quad 
\pi(a)=a\quad\forall a\in A, 
\]
where $A$ is identified with subalgebras of the crossed products. 
For $\gamma\in\Aut(C^*(A,\phi^u_{(1,0)}))$, one can check 
\[
\zeta_{\phi_{(1,0)}}(\tilde\eta_0(\pi^{-1}\circ\gamma\circ\pi))
=\zeta_{\phi^u_{(1,0)}}(\tilde\eta_0(\gamma))
\in\OrderExt(K_0(A),K_0(A),D_A), 
\]
where $\tilde\eta_0$ in the left hand side is defined 
for $C^*(A,\phi_{(1,0)})$ and 
$\tilde\eta_0$ in the right hand side is defined 
for $C^*(A,\phi^u_{(1,0)})$. 
We also have 
\[
(\pi^{-1}\circ\tilde\phi^u_{(0,1)}\circ\pi)(a)
=\pi^{-1}(\tilde\phi^u_{(0,1)}(a))=\tilde\phi^u_{(0,1)}(a)
=(\Ad u_{(0,1)}\circ\tilde\phi_{(0,1)})(a)\quad\forall a\in A
\]
and 
\begin{align*}
(\pi^{-1}\circ\tilde\phi^u_{(0,1)}\circ\pi)(\lambda_{\phi_{(1,0)}})
&=(\pi^{-1}\circ\tilde\phi^u_{(0,1)})(u_{(1,0)}^*\lambda_{\phi^u_{(1,0)}}) \\
&=\pi^{-1}(\tilde\phi^u_{(0,1)}(u_{(1,0)}^*)\lambda_{\phi^u_{(1,0)}}) \\
&=\phi^u_{(0,1)}(u_{(1,0)}^*)u_{(1,0)}\lambda_{\phi_{(1,0)}} \\
&=u_{(0,1)}\phi_{(0,1)}(u_{(1,0)}^*)u_{(0,1)}^*
u_{(1,0)}\lambda_{\phi_{(1,0)}} \\
&=u_{(0,1)}\phi_{(1,0)}(u_{(0,1)}^*)u_{(1,0)}^*
u_{(1,0)}\lambda_{\phi_{(1,0)}} \\
&=u_{(0,1)}\phi_{(1,0)}(u_{(0,1)}^*)\lambda_{\phi_{(1,0)}} \\
&=u_{(0,1)}\lambda_{\phi_{(1,0)}}u_{(0,1)}^* \\
&=(\Ad u_{(0,1)}\circ\tilde\phi_{(0,1)})(\lambda_{\phi_{(1,0)}}). 
\end{align*}
Thus 
$\pi^{-1}\circ\tilde\phi^u_{(0,1)}\circ\pi
=\Ad u_{(0,1)}\circ\tilde\phi_{(0,1)}$. 
Since inner automorphisms are contained in the kernel of $\tilde\eta_0$, 
we obtain 
\begin{align*}
\zeta_{\phi_{(1,0)}}(\tilde\eta_0(\tilde\phi_{(0,1)}))
&=\zeta_{\phi_{(1,0)}}(\tilde\eta_0(\Ad u_{(0,1)}\circ\tilde\phi_{(0,1)})) \\
&=\zeta_{\phi_{(1,0)}}(\tilde\eta_0
(\pi^{-1}\circ\tilde\phi^u_{(0,1)}\circ\pi)) \\
&=\zeta_{\phi^u_{(1,0)}}(\tilde\eta_0(\tilde\phi^u_{(0,1)})), 
\end{align*}
which completes the proof. 
\end{proof}

\begin{thm}\label{Z2cc}
Let $\phi,\psi:\Z^2\curvearrowright A$ be 
uniformly outer, locally $KK$-trivial $\Z^2$-actions 
on a unital simple AF algebra $A$ with unique trace. 
The following are equivalent. 
\begin{enumerate}
\item $[\phi]=[\psi]$. 
\item $\phi$ and $\psi$ are $KK$-trivially cocycle conjugate. 
\end{enumerate}
\end{thm}
\begin{proof}
(2)$\Rightarrow$(1) was shown in the proposition above 
without assuming that $A$ has a unique trace. 
Let us consider the other implication (1)$\Rightarrow$(2). 
By Theorem \ref{ZRohlintype} and Theorem \ref{Zcc}, 
we may assume that there exists a unitary $u\in A$ 
such that $\psi_{(1,0)}=\Ad u\circ\phi_{(1,0)}$. 
By Theorem \ref{ZRohlintype} and Remark \ref{modelaction} 
(or Remark \ref{Lin}), 
the crossed product $C^*$-algebra $C^*(A,\phi_{(1,0)})$ is 
a unital simple AT algebra with real rank zero. 

Clearly $\phi_{(0,1)}$ extends to 
$\tilde\phi_{(0,1)}\in\Aut(C^*(A,\phi_{(1,0)}))$ by 
\[
\tilde\phi_{(0,1)}(a)=a\quad\forall a\in A\quad\text{and}\quad 
\tilde\phi_{(0,1)}(\lambda_{\phi_{(1,0)}})=\lambda_{\phi_{(1,0)}}. 
\]
Since
\[
\psi_{(0,1)}\circ\phi_{(1,0)}
=\Ad(\psi_{(0,1)}(u^*)u)\circ\phi_{(1,0)}\circ\psi_{(0,1)}, 
\]
we can extend $\psi_{(0,1)}$ to $\omega\in\Aut(C^*(A,\phi_{(1,0)}))$ by 
\[
\omega(a)=a\quad\forall a\in A\quad\text{and}\quad 
\omega(\lambda_{\phi_{(1,0)}})=\psi_{(0,1)}(u^*)u\lambda_{\phi_{(1,0)}}. 
\]
In order to apply Theorem \ref{equivEK} 
to $\tilde\phi_{(0,1)}$ and $\omega$, 
we would like to check that 
these automorphisms are asymptotically unitarily equivalent. 
There exists an isomorphism $\pi:C^*(A,\phi_{(1,0)})\to C^*(A,\psi_{(1,0)})$ 
defined by 
\[
\pi(a)=a\quad\forall a\in A\quad\text{and}\quad 
\pi(\lambda_{\phi_{(1,0)}})=u^*\lambda_{\psi_{(1,0)}}. 
\]
As mentioned in the proof of Proposition \ref{invariance}, 
for any $\gamma\in\Aut(C^*(A,\phi_{(1,0)}))$, one has 
\[
\zeta_{\phi_{(1,0)}}(\tilde\eta_0(\gamma))
=\zeta_{\psi_{(1,0)}}(\tilde\eta_0(\pi\circ\gamma\circ\pi^{-1})). 
\]
Moreover it is easy to see that 
$\pi\circ\omega\circ\pi^{-1}$ is equal to $\tilde\psi_{(0,1)}$, 
which is defined by 
\[
\tilde\psi_{(0,1)}(a)=a\quad\forall a\in A\quad\text{and}\quad 
\tilde\psi_{(0,1)}(\lambda_{\psi_{(1,0)}})=\lambda_{\psi_{(1,0)}}. 
\]
It follows that 
\begin{align*}
\zeta_{\phi_{(1,0)}}(\tilde\eta_0(\omega))
&=\zeta_{\psi_{(1,0)}}(\tilde\eta_0(\pi\circ\omega\circ\pi^{-1})) \\
&=\zeta_{\psi_{(1,0)}}(\tilde\eta_0(\tilde\psi_{(0,1)}))
=[\psi]
=[\phi]
=\zeta_{\phi_{(1,0)}}(\tilde\eta_0(\tilde\phi_{(0,1)})), 
\end{align*}
and so $\tilde\eta_0(\omega)=\tilde\eta_0(\tilde\phi_{(0,1)})$. 
By Lemma \ref{eta1=0}, 
$\eta_1(\omega)=\eta_1(\tilde\phi_{(0,1)})=0$. 
Therefore, by Theorem \ref{KishiKum}, $\tilde\phi_{(0,1)}$ and $\omega$ are 
asymptotically unitarily equivalent. 

Then, Theorem \ref{equivEK} applies and yields 
an approximately inner automorphism $\mu\in\Aut_\T(C^*(A,\phi_{(1,0)}))$ 
and a unitary $v\in A$ such that 
$\mu|A$ is in $\aInn(A)$ and 
\begin{equation}
\mu\circ\omega\circ\mu^{-1}=\Ad v\circ\tilde\phi_{(0,1)}. 
\label{cat}
\end{equation}
By restricting this equality to $A$, we get 
\begin{equation}
(\mu|A)\circ\psi_{(0,1)}\circ(\mu|A)^{-1}
=\Ad v\circ\phi_{(0,1)}. 
\label{cow}
\end{equation}
Let $z\in A$ be the unitary satisfying 
$\mu(\lambda_{\phi_{(1,0)}})=z\lambda_{\phi_{(1,0)}}$. 
Then 
\begin{equation}
(\mu|A)\circ\psi_{(1,0)}\circ(\mu|A)^{-1}
=(\mu|A)\circ\Ad u\circ\phi_{(1,0)}\circ(\mu|A)^{-1}
=\Ad\mu(u)z\circ\phi_{(1,0)}. 
\label{tiger}
\end{equation}
From \eqref{cat}, one can see that 
\begin{align*}
(\mu\circ\omega\circ\mu^{-1})(\lambda_{\phi_{(1,0)}})
&=(\mu\circ\omega)(\mu^{-1}(z^*)\lambda_{\phi_{(1,0)}}) \\
&=\mu(\psi_{(0,1)}(\mu^{-1}(z^*))\psi_{(0,1)}(u^*)u\lambda_{\phi_{(1,0)}}) \\
&=(\Ad v\circ\phi_{(0,1)})(z^*\mu(u^*))\mu(u)z\lambda_{\phi_{(1,0)}} \\
&=v\phi_{(0,1)}(z^*\mu(u^*))v^*\mu(u)z\lambda_{\phi_{(1,0)}}
\end{align*}
is equal to 
\begin{align*}
(\Ad v\circ\tilde\phi_{(0,1)})(\lambda_{\phi_{(1,0)}})
&=v\lambda_{\phi_{(1,0)}}v^* \\
&=v\phi_{(1,0)}(v^*)\lambda_{\phi_{(1,0)}}. 
\end{align*}
Hence one obtains 
\begin{equation}
v\phi_{(0,1)}(\mu(u)z)=\mu(u)z\phi_{(1,0)}(v). 
\label{rabbit}
\end{equation}
It follows from \eqref{cow}, \eqref{tiger}, \eqref{rabbit} that 
$\psi$ and $\phi$ are $KK$-trivially cocycle conjugate. 
\end{proof}

\begin{rem}
We do not know 
the precise range of our invariant which takes its values in $\OrderExt$. 
At least, the following observation shows that 
the range does not exhaust $\OrderExt$. 
Let $\phi:\Z^2\curvearrowright A$ be 
a locally $KK$-trivial and uniformly outer $\Z^2$-action 
on a unital simple AF algebra. 
Suppose that 
$(\xi,R)$ is a representative of $[\phi]\in\OrderExt(K_0(A),K_0(A),D_A)$. 
Since 
\[
\begin{CD}
\xi:\quad 0 @>>> K_0(A) @>\iota>> E_\xi @>q>> K_0(A) @>>> 0
\end{CD}
\]
is exact and $R:E_\xi\to\Aff(T(A))$ satisfies $R\circ\iota=D_A$, 
there exists a homomorphism $R_0:K_0(A)\to\Aff(K_0(A))/\Ima D_A$ such that 
$R_0(q(x))=R(x)+D_A(K_0(A))$ for any $x\in E_\xi$. 
It is easy to see $R_0([1_A])=0$, 
because the implementing unitary $\lambda_{\phi_{(1,0)}}$ 
is fixed by $\tilde\phi_{(0,1)}$. 
Thus, $[\phi]$ belongs to the subgroup 
\[
\{[(\xi,R)]\in\OrderExt(K_0(A),K_0(A),D_A)\mid R_0([1_A])=0\}. 
\]
When $A$ is a UHF algebra, one can see that 
this subgroup coincides 
with the range of the invariant introduced in \cite{KM}. 
Therefore, 
Theorem \ref{Z2cc} yields a new proof of \cite[Theorem 6.5]{KM}. 
\end{rem}

\end{document}